\documentclass[12pt]{iopart}

\usepackage{iopamsm}
 \usepackage[noadjust]{cite}
 \usepackage{color}
 \usepackage{curves}
 \usepackage{enumerate}
 \usepackage{epic}
 \usepackage{epsfig}
 \usepackage{float}
 \usepackage{graphics}
 \usepackage{graphicx}
 \usepackage{latexsym}
 \usepackage{multicol}
 \usepackage{subfigure}
\usepackage{theorem}
 \newtheorem{thm}{Theorem}[section]
 \newtheorem{hyp}{Hypothesis}

 \newtheorem{prop}[thm]{Proposition}
 \theorembodyfont{\upshape}
 
\numberwithin{equation}{section}
 \DeclareMathOperator{\A}{A}

 \DeclareMathOperator{\ccc}{c.c.}

 \DeclareMathOperator{\real}{Re}

 \DeclareMathOperator{\GL}{GL}
 
 \newcommand{\proof}{\textbf{Proof:\quad}}

 \newcommand{\GC}{\mathbb{C}}

 \newcommand{\GM}{\mathbb{M}}
 \newcommand{\GN}{\mathbb{N}}

 \newcommand{\GR}{\mathbb{R}}
 \newcommand{\GS}{\mathbb{S}}
 \newcommand{\GSE}{\mathbb{SE}}
 \newcommand{\GSO}{\mathbb{SO}}
 
 \newcommand{\GZ}{\mathbb{Z}}

 \newcommand{\omr}{\omega_{\tiny\mbox{rot}}}
 \newcommand{\ov}{\overline}

 \newcommand{\sdp}{\dot{+}}

 \newcommand{\newl}{\newline\newline}
  \DeclareRobustCommand{\qed}{%
  \ifmmode \mathqed
  \else
    \leavevmode\unskip\penalty9999 \hbox{}\nobreak\hfill
    \quad\hbox{$\Box$\normalsize}%
  \fi
}

\renewcommand{\footnoterule}{%
\kern -3pt 
\hrule height 0.6pt width 0.4\columnwidth
\kern 2.6pt 
}
\begin{document}

\title[Spiral anchoring in anisotropic media with multiple inhomogeneities]{Spiral anchoring in anisotropic media with multiple inhomogeneities}

\author{P Boily}

\address{Department of Mathematics and Statistics, University of Ottawa,
Ottawa K1N 6N5, Canada} \ead{pboily@uottawa.ca}
\begin{abstract}
Various PDE models have been suggested in order to explain and
predict the dynamics of spiral waves in excitable media. In two
landmark papers, Barkley noticed that some of the behaviour could be
explained by the inherent Euclidean symmetry of these models.
LeBlanc and Wulff then introduced forced Euclidean symmetry-breaking
(FESB) to the analysis, in the form of individual translational
symmetry-breaking (TSB) perturbations and rotational
sym\-metry-breaking (RSB) perturbations; in either case, it is shown
that spiral anchoring is a direct consequence of the FESB.
\par
In this article, we provide a characterization of spiral anchoring
when two perturbations, a TSB term and a RSB term, are combined,
where the TSB term is centered at the origin and the RSB term
preserves rotations by multiples of $\frac{2\pi}{\jmath^*}$, where
$\jmath^*\geq 1$ is an integer. When $\jmath^*
>1$ (such as in a modified bidomain model), it is shown that spirals anchor at the origin, but when
$\jmath^* =1$ (such as in a planar reaction-diffusion-advection
system), spirals generically anchor away from the origin.
\end{abstract}


\maketitle

\section{Introduction}
Scientific literature has recorded the occurrences of spiral waves
in excitable media, from the Belousov-Zhabotinsky chemical reaction
to the electrical potential in cardiac tissue
\cite{LOPS,GZM,ZM,MPMPV,WR,YP,Detal,R1,Wetal,J,MAK,B1,B2,BKT}. This
last instance provides an inspiration for the study of spiral wave,
as they are believed to be a precursor to several fatal types of
ventricular tachycardia and/or fibrillation. \cite{W,Wetal,KS}.

After Barkley realized that the resonant growth transition from
rigid rotations to quasi-periodic meandering and drifting could be
explained solely through the inherent Euclidean symmetries of
excitable media (which were incorporated in the various models)
\cite{B1,B2,BK}, the theory of equivariant dynamical systems
(together with the general center manifold reduction theorem of
Sandstede, Scheel and Wulff \cite{SSW1,SSW2,SSW3,SSW4,FSSW}) has
been used on numerous occasions to explain and predict spiral wave
dynamics \cite{GLM,Wulff,GLM2,LW,LeB,C1,BLM,Bo2,Byeah}.

One of the advantages of this approach is that it provides
universal, model-independent explanations of observed dynamics and
bifurcations of spiral waves that occur in \textit{a priori}
different situations; consider for instance the aforementioned Hopf
bifurcation from rigid rotation to quasi-periodic meandering and
drifting (observed in numerical simulations \cite{BKT} and in actual
chemical reactions \cite{LOPS}) or the anchoring/repelling of spiral
waves on/from a site of inhomogeneity (seen in numerical
integrations of an Oregonator system \cite{MPMPV}, in
photo-sensitive chemical reactions \cite{ZM} and in cardiac tissue
\cite{Detal}). Using a model-independent approach based on forced
symmetry-breaking, it has been shown that anchoring/repelling of
rotating waves generically takes place in systems in which the
translation symmetry of $\mathbb{S}\mathbb{E}(2)$ (the group of
translation and rotation symmetries) is broken by a small
perturbation \cite{LW}, although the center of anchoring/repelling
may not be an inhomogeneity when more than one such perturbations
are present \cite{BLM}. Similarly, some dynamics of spiral waves
observed in anisotropic media ({\em e.g.} phase-locking and/or
linear drifting of meandering spiral waves) have been shown to be
generic consequences of rotational symmetry-breaking
\cite{R1,R2,LeB}. \par Many of the phenomena in which spiral waves
are observed experimentally are modeled by reaction-diffusion
systems
\begin{equation}
\frac{\partial u}{\partial t}=D\cdot\nabla^2\,u+f(u)
\label{basicrdpde}
\end{equation}
where $u$ is a $k$-vector valued function of time and
two-dimensional space, $D$ is a matrix of diffusion coefficients and
$f:\mathbb{R}^k\longrightarrow \mathbb{R}^k$ is a smooth reaction
term. Their use as models is justified as Scheel has proved that
systems of this form admit have time-periodic, rigidly rotating
spiral wave solutions \cite{S}. Such models are
$\GSE(2)-$equivariant: indeed, (\ref{basicrdpde}) is invariant under
the transformations
\begin{equation}
u(t,x)\longmapsto
u(t,x_1\cos\,\theta-x_2\sin\,\theta+p_1,x_1\sin\,\theta+
x_2\cos\,\theta+p_2), \label{uaction}
\end{equation}
where
$(\theta,p_1,p_2)\in\mathbb{S}^1\times\mathbb{R}^2\simeq\GSE(2)$ and
$x\in \GR^2$ \cite{Wulff,DMcK}. But inhomogeneous media is not
perfectly Euclidean: as a result, (\ref{basicrdpde}) does not
provide an appropriate frame for the study of anisotropic media with
inhomogeneities.
\par
Consider rather the following generalization of the bidomain
equations describing the electrical properties of anisotropic
cardiac tissue
\begin{align}\label{thebidomain}
\begin{split}
u_t&=\textstyle{\frac{1}{\varsigma}(u-\frac{u^3}{3}-v)+\nabla^2 u+\frac{\alpha \varepsilon}{1+\alpha (1-\varepsilon)}\Psi_{x_1x_1}}\\
v_t&=\varsigma(u+\beta-\gamma v),\\
\nabla^2 &\Psi+\varepsilon g(\alpha,\varepsilon)\Psi_{x_2x_2}=
\varepsilon h(\alpha,\varepsilon) u_{x_2x_2},
\end{split}
\end{align}
where $u$ is a transmembrane potential, $v$ controls the recovery of
the action potential, $\Psi$ is an auxiliary potential (without
obvious physical interpretation), $x_1$ is the preferred direction
in physical space in which tissue fibers align, $\varepsilon$~is a
measure of that preference, $g$ and $h$ are appropriate model
functions, and $\alpha,\varsigma,\beta$ and $\gamma$ are model
parameters \cite{R2,LR,BEL,Byeah}.\par If the tissue has equal
anisotropy ratios (\textsl{i.e.} $\varepsilon=0$),
(\ref{thebidomain}) decouples into the FitzHugh-Nagumo equations for
$u$ and $v$, and Poisson's equation for $\Psi$ \cite{R2}. Ignoring
the boundary, system (\ref{thebidomain}) is $\GSE(2)-$equivariant
when $\varepsilon=0$, while it is only $\GZ_2\sdp \GR^2-$equivariant
otherwise (see \cite[section 2.2]{LR} for details). \par Cardiac
tissue is littered with inhomogeneities. In order to make the
analysis more tractable, let us make the modeling assumption that
the inhomogeneities consist of a finite number of independent
``sources'' which are localized near distinct sites $\zeta_1,\ldots,
\zeta_n$ in the plane (see \cite{BLM} for a similar hypothesis).
This situation could be modeled by the (slightly) perturbed bidomain
equations
\begin{align}\label{thebidomain2}
\begin{split}
u_t&=\frac{1}{\varsigma}(u-\frac{u^3}{3}-v)+\nabla^2 u+\frac{\alpha
\varepsilon}{1+\alpha (1-\varepsilon)}\Psi_{x_1x_1}+
\sum_{j=1}^n\mu_j g_j^{u}(\|x-\zeta_j\|^2,\mu)\\
v_t&=\varsigma(u+\beta-\gamma v)+\sum_{j=1}^n\mu_j g_j^{v}(\|x-\zeta_j\|^2,\mu), \\
\nabla^2& \Psi+\varepsilon g(\alpha,\varepsilon)\Psi_{x_2x_2}=
\varepsilon h(\alpha,\varepsilon) u_{x_2x_2},
\end{split}
\end{align} where $\mu=(\mu_1,\ldots,\mu_n)\in \GR^n$ is a small parameter and $g_j^{u,v}$ are smooth functions, uniformly bounded in their variables.
\par The symmetry-breaking induced by anisotropy and the presence of
inhomogeneities at the origin then yields an equivariant structure
on the semi-flow $\Phi_{t,\varepsilon,\mu}$ of (\ref{thebidomain2}),
the existence and uniqueness of which is shown in \cite{FS}. In
particular, $\Phi_{t,\varepsilon,\mu}$
\begin{enumerate}[(E1)]
\item is $\GSE(2)-$equivariant when $(\varepsilon,\mu)=0$;
\item is $\GZ_{2}\sdp\GR^2-$equivariant when $\varepsilon\neq 0$ is small and $\mu=0$;
\item preserves rotations around $\zeta_{i^*}$ (but generically not translations) when
$\varepsilon=0$ and all the $\mu_i$ are zero except $\mu_{i*}$, and
\item is (generically) trivially equivariant when $(\varepsilon,\mu)$ is small and
generic.
\end{enumerate}
\par When $\varepsilon=0$ and $n=1$, LeBlanc and Wulff have shown that rotating waves generically anchor at (or are repelled by)
the inhomogeneity \cite{LW}. When $\varepsilon=0$ and $n>1$, Boily,
LeBlanc and Matsui have shown that anchoring remains generic in
parameter wedges, but that the anchoring/repelling is generically
centered away from the inhomogeneities \cite{BLM}. When $\mu=0$,
LeBlanc has shown that rotating waves have $\GZ_2-$spatial symmetry
\cite{LeB}.
\par There is no reason to believe that combining anisotropy with
multiple inhomogeneities will change these results qualitatively: to
wit, if a spiral anchors in non-anisotropic media, it should also do
so in anisotropic media, but with $\GZ_2-$spatial symmetry. However,
the proofs of \cite{LW,BLM,LeB} cannot be pasted together to obtain
this ``obvious'' result; some mathematical difficulties had to be
overcome.
\newl The goal of this paper, then, is to provide a detailed
analysis of abstract dynamical systems which share properties (E1),
(E3) and (E4), as well as some technical conditions which will be
specified as we proceed, with the following slight modification:
\begin{enumerate}[(E2')]
\item the semi-flow is $\GZ_{\jmath^*}\sdp\GR^2-$equivariant when $\varepsilon\neq 0$ is small and
$\mu=0$, where $\jmath^*$ is some positive integer.
\end{enumerate}
As noted, for the bidomain model we have $\jmath^*=2$, and for a
reaction-diffusion-advection system, we have $\jmath^*=1$. The value
of $\jmath^*$ will play an important role in the subsequent
analysis. The paper is organized as follows. In the second section,
we derive the center bundle equations of the semi-flow of an
appropriate system near a hyperbolic rotating wave. We state and
prove our main results for $n=1$ in the third section: as long as
$\jmath^*>1$, spirals anchor qualitatively as they would in the
absence of anisotropy. But in the case $\jmath^*=1$, an unexpected
result occurs: in a reaction-diffusion-advection system with one
inhomogeneity, the tug of war between the advection current and the
inhomogeneity ends in a tie, and spirals anchor at some point away
from the inhomogeneity. To the best of our knowledge, this has yet
to be observed either in simulations or in the lab. In the fourth
section, we look at spiral anchoring in the presence of the most
general form of Euclidean symmetry-breaking. Finally, we perform a
simple qualitative numerical experiment demonstrating the validity
of our results.


\section{Reduction to the Center Bundle Equations}

Set $1\leq \jmath^*\in \GN$. Let $X$ be a Banach space, ${\mathcal
U}\subset\GR\times\mathbb{R}^n$ a neighborhood of the origin, and
$\Phi_{t,\varepsilon,\mu}$ be a smoothly parameterized family (parameterized by
$(\varepsilon,\mu)\in {\mathcal U}$) of smooth local semi-flows on $X$.

Let $\mathbb{S}\mathbb{E}(2)=\mathbb{C}\dot{+}\mathbb{S}\mathbb{O}(2)$ denote
the group of all planar translations and rotations, and let
\begin{equation}
a:\mathbb{S}\mathbb{E}(2)\longrightarrow \GL(X) \label{a_action}
\end{equation}
be a faithful and isometric representation of $\mathbb{S}\mathbb{E}(2)$ in the
space of bounded, invertible linear operators on $X$.  For example, if $X$ is a
space of functions with planar domain, a typical $\mathbb{S}\mathbb{E}(2)$
action (such as (\ref{uaction}) in the preceding section) is given by

$$(a(\gamma)u)(x)=u(\gamma^{-1}(x)),\,\,\,\,\gamma\in\mathbb{S}\mathbb{E}(2).$$
\par We will parameterize $\mathbb{S}\mathbb{E}(2)$ as follows:
$\mathbb{S}\mathbb{E}(2)\cong\,\mathbb{C}\times\mathbb{S}^1$, with
multiplication given by $(p_1,\varphi_1)\cdot (p_2,\varphi_2)=( e^{
i \varphi_1}p_2+p_1,\varphi_1+\varphi_2)$,
$\forall\,(p_1,\varphi_1),\,(p_2,\varphi_2)\in\,\mathbb{C}\times\mathbb{S}^1$.
For fixed $\xi\in\mathbb{C}$, define the following subgroups of
$\mathbb{S}\mathbb{E}(2)$:
$$\GC\sdp\mathbb{Z}_{\jmath^*}
=\left\{\,\left(s,k\frac{2\pi}{\jmath^*}\right)\,\,|\,\, k\in \GZ,\
s\in \GC\,\right\}\mbox{
 and  }\mathbb{S}\mathbb{O}(2)_{\xi}=\{\,(\xi,0)\cdot
(0,\theta)\cdot (-\xi,0)\,\,|\,\,\theta \in\mathbb{S}^1\,\},$$ the
latter of which is isomorphic to $\mathbb{S}\mathbb{O}(2)$ and
represents rotations about the point $\xi$.
  We will assume the following symmetry conditions on
the family $\Phi_{t,\varepsilon,\mu}$ of semi-flows.
\begin{hyp} \label{hyp1}
There exists $n$ distinct points $\xi_1,\ldots,\xi_n$ in $\mathbb{C}$ such that
if $e_j$ denotes the $j^{\mbox{\footnotesize th}}$ vector of the canonical
basis in $\mathbb{R}^n$, then $\forall\,u\in\,X,\ \varepsilon\neq 0,\
\alpha\neq 0,\  t>0,$
\begin{align*}
\Phi_{t,\varepsilon,0}(a(\gamma)u)&=a(\gamma)\Phi_{t,\varepsilon,0}(u) \iff
\gamma\in\,\GC\sdp \mathbb{Z}_{\jmath^*},\\
\Phi_{t,0,\alpha e_j}(a(\gamma)u)&=a(\gamma)\Phi_{t,0,\alpha
e_j}(u) \iff \gamma\in\,\mathbb{S}\mathbb{O}(2)_{\xi_j}, \quad\mbox{and}\\
\Phi_{t,0,0}(a(\gamma)u)&=a(\gamma)\Phi_{t,0,0}(u), \quad \forall\,
\gamma\in\,\mathbb{S}\mathbb{E}(2).
\end{align*}
\end{hyp}
Hypothesis \ref{hyp1} basically states that
$\Phi_{t,\varepsilon,\mu}$ satisfies properties (E1), (E2'), (E3),
(E4). We are interested in the effects of the forced
symmetry-breaking on normally hyperbolic rotating waves. Therefore,
we will assume the following hypothesis.
\begin{hyp}
There exists $u^*\in X$ and $\Omega^*$ in the Lie algebra of
$\mathbb{S}\mathbb{E}(2)$ such that $e^{\Omega^*t}$ is a rotation and
$\Phi_{t,0}(u^*)=a(e^{\Omega^*t})u^*$ for all $t$. Moreover, the set
$\{\,\lambda\in\mathbb{C}\,\,|\,\,|\lambda|\geq 1\,\}$ is a spectral set for
the linearization $a(e^{-\Omega^*})D\Phi_{1,0}(u^*)$ with projection $P_*$ such
that the generalized eigenspace $\mbox{\rm range}(P_*)$ is three dimensional.
\label{hyp2}
\end{hyp}
In order to simplify the analysis, we only consider one-armed spiral
waves, \textit{i.e.} the isotropy subgroup of $u^*$ in hypothesis
\ref{hyp2} is trivial. It should be noted that hypotheses 1 and 2
hold for a large variety of spirals (such as decaying spirals), but
not for all spirals (including Archimedean spirals)
\cite{S}.\footnote{Even when hypothesis 1 fails, finite-dimensional
center-bundle equations which share the symmetries of the underlying
abstract dynamical systems have a definite predictive value
\cite{B1,B2,LOPS,LW,LeB}.}
\par
From now on, we assume $\Phi_{t,\varepsilon,\mu}$ is a semi-flow
that satisfy both hypotheses, as well as all other hypotheses
required in order for the center manifold theorems of
\cite{SSW1,SSW2,SSW3,SSW4} to hold. Then, for $(\varepsilon,\mu)$
near the origin in $\GR\times \mathbb{R}^n$, the essential dynamics
of the semi-flow $\Phi_{t,\varepsilon,\mu}$ near the rotating wave
reduces to the following ordinary differential equations on the
bundle $\mathbb{C}\times\mathbb{S}^1$ (see \cite{Byeah} for more
details):
\begin{align}\label{basic10}
\begin{split}
\dot{p}&= e^{ i \varphi}\big[\nu+J^p(p,\ov{p},\varphi,\varepsilon,\mu)\big]\\
\dot{\varphi}&=\omr+J^{\varphi}(p,\ov{p},\varphi,\varepsilon,\mu),
\end{split}
\end{align}
where $\nu$ is a complex constant, $0\neq\omr$ is a real constant,
$J^p(p,\ov{p},\varphi,0,0)\equiv 0$ and
$J^{\varphi}(p,\ov{p},\varphi,0,0)\equiv 0$. Furthermore, the functions $J^p$
and $J^{\varphi}$ are smooth and uniformly bounded in $p$.  If
$(\varepsilon,\mu)$ is near the origin, we can re-scale time along orbits of
(\ref{basic10}), perform a simple computation and apply Taylor's theorem to get
the following.
\begin{prop}
The symmetry conditions in hypothesis $\ref{hyp1}$ imply that the equations
$(\ref{basic10})$ have the general form
\begin{align}\label{basic12}
\begin{split}
\dot{p}&= e^{ i \varphi(t)}\left[v+\varepsilon
G(\varphi(t),\varepsilon)+\sum_{j=1}^n\mu_j H_j((p-\xi_j) e^{- i
\varphi(t)},\ov{(p-\xi_j)} e^{ i \varphi(t)},\mu_j)\right]
\end{split}
\end{align}
where, without loss of generality, $\varphi(t)=t$, $v\in \GC$,
$\mu=(\mu_1,\ldots,\mu_n)$, and the functions $G,H_j$ are smooth,
periodic in $\varphi$ and uniformly bounded in $p$, and $G$ is
$2\pi/\jmath^*-$periodic in $\varphi$ for some positive integer
$\jmath^*$.\footnote{Strictly speaking, the center bundle equations
take the form
\begin{align}\label{whatsthedeal}\dot{p}=\mathcal{F}+e^{it}\sum_{j\neq
k}
\lambda_j\lambda_k\mathcal{H}_{j,k}(p,\ov{p},\xi_j,\xi_k,t,\lambda),\end{align}
where $\mathcal{F}$ is the right hand side of (\ref{basic12}),
$\lambda=(\lambda_0,\lambda_1,\ldots,\lambda_n)=(\varepsilon,\mu_1,\ldots,\mu_n)\in
\GR^{n+1}$, and the functions $\mathcal{H}_{j,k}$ are suitably
smooth and $2\pi-$periodic in $t$. However, the analysis in what
follows is independent of the ``mixed'' perturbation terms
$\mathcal{H}_{j,k}$, when $\|\lambda\|$ is small enough: in
particular, it is then generically the case that if
(\ref{whatsthedeal}) has a $2\pi-$periodic solution, so does
(\ref{basic12}), and vice-versa, and these solutions share the same
stability. The argument depends on the Taylor expansion of an the
time$-2\pi$ map (see \cite{BLM} for an example of a similar analysis
that includes the mixed perturbation terms).}
\end{prop}
As in \cite{BLM}, a $2\pi-$periodic solution $p_{\varepsilon,\mu}$ of
$(\ref{basic12})$ is called a \textit{perturbed rotating wave} of
$(\ref{basic12})$. Define the average value
\begin{align}[p_{\varepsilon,\mu}]_{\A}&=\frac{1}{2\pi}\int_{0}^{2\pi}\!\!\!\!p_{\varepsilon,\mu}(t)\,  d t.\label{ca}\end{align}
If the Floquet multipliers of $p_{\varepsilon,\mu}$ all lie within (resp.
outside) the unit circle, we shall say that $[p_{\varepsilon,\mu}]_{\A}$ is the
\textit{anchoring} (resp. \textit{repelling}, or \textit{unstable anchoring})
\textit{center} of~$p_{\varepsilon,\mu}$.\par In the following section, we
 perform an analysis of the anchoring of perturbed rotating waves of
(\ref{basic12}) for $n=1$ and $\xi_1=0$ and parameter values
$(\varepsilon,\mu)$ near the origin in $\GR^2$. As noted, we will
tackle the case $n>1$ in section~\ref{general}.


\section{Analysis of the Center Bundle Equations $(n=1,\xi_1=0)$}

Equations (\ref{basic12}) represent the dynamics near a normally
hyperbolic rotating wave for a parameterized family
$\Phi_{t,\varepsilon,\mu}$ of semi-flows satisfying the
forced-symmetry breaking conditions in hypothesis \ref{hyp1}. We
start with a brief review of spiral anchoring in the two cases
$\varepsilon=0$ and $\mu_1=0$, which were studied in detail in
\cite{LW,LeB,BLM}, and then present the result of combining both
types of perturbations.
\subsection{The Case $\varepsilon=0$}\label{thecasee=0}
By writing $w=p  e^{- i t}\!\!+i v$, this system becomes
\begin{align}
\label{basiceqs5} \dot{w}&={\displaystyle -i
w+\mu_1{\widetilde{H}}(w,\overline{w},\mu_1)}
\end{align}
where $\widetilde{H}(w,\overline{w},\mu_1)= H_1(w-i
v,\overline{w}+i\overline{v},0,\mu_1)$. The following theorem is
proved in \cite{LW}.
\begin{thm}
Let $\alpha=\mbox{\rm Re}(D_1\widetilde{H}(0,0,0))$, where
$\tilde{H}$ is as in $(\ref{basiceqs5})$.  If $\alpha\neq 0$, then
for all $\mu_1\neq 0$ small enough,
 $(\ref{basic12})$ has a hyperbolic rotating wave
\begin{equation}
p(t)=\left(-i v+ O(\mu_1)\right)  e^{ i  t},\quad\varphi(t)=t.
\label{anchored_sol}
\end{equation} The origin $[p]_{\A}=0$ is an anchoring center if $a\mu_1<0$; it is a
repelling center if $a\mu_1>0$.
\end{thm}
In the case where the semi-flow $\Phi_{t,0,\mu_1}$ is generated by a
system of planar reaction-diffusion partial differential equations,
the solution (\ref{anchored_sol}) represents a wave which is rigidly
and uniformly rotating around the origin in the plane. In the case
where $\alpha\mu_1<0$, the rotating wave is locally asymptotically
stable. When $\alpha\mu_1>0$, the rotating wave is unstable (see
\cite{MPMPV} for an experimental characterization of this phenomenon
in an Oregonator model).
\subsection{The Case $\mu_1=0$} Since $G$ is $2\pi/\jmath^*-$periodic in $\varphi$ (or
in $t$, equivalently), it can be written as the uniformly convergent Fourier
series
\begin{align}\label{FourierG1}
G(t,\varepsilon)=\sum_{m\in {\mathbb Z}}\, g_m(\varepsilon) e^{ i
m\jmath^*t}.
\end{align}
The following theorem is proved in \cite{LeB}.
\begin{thm}
Let $\jmath^*>1$, or $\jmath^*=1$ and $g_{-1}(\varepsilon)\equiv 0$.
Then, for all $\varepsilon\neq 0$ sufficiently small, the solutions
of $(\ref{basic12})$ are
 $2\pi/\jmath^*-$periodic in time with discrete $\GZ_{{\jmath^*}}-$symmetry.
\end{thm}
In the case where the semi-flow $\Phi_{t,\varepsilon,0}$ is generated by a
system of planar reaction-diffusion partial differential equations, the
solutions of (\ref{basic12}) represents discrete rotating waves in the physical
space.\newl How do these two cases interact when they are combined? We shall
see that the answer depends greatly on the nature of $\jmath^*$.
\subsection{The Case $\jmath^*=1$} \label{thecasel=1}
Let $F_G:\GR\times\GR^2\to \GC$ be defined by
\begin{equation}\label{FG1}
F_G(t,\varepsilon)= e^{ i t}\Big[- i v+\varepsilon \sum_{m\neq
-1}\,\frac{g_m(\varepsilon) e^{ i m t}}{ i(m+1)}\Big],
\end{equation} where the $g_m$ are the Fourier coefficients of $G$ found in
(\ref{FourierG1}). Differentiating (\ref{FG1}) yields
\begin{align}\label{FGdot1}
\dot{F}_G(t,\varepsilon)&= e^{ i t}\big[v+\varepsilon
G(t,\varepsilon)-\varepsilon g_{-1}(\varepsilon)\big].
\end{align}
Generically, on the center bundle $\GC\times\GS^1$, the path
$(F_G(t,\varepsilon),t)$ represents a discrete rotating wave around
the origin with trivial spatio-temporal symmetry. Set $z=p-F_G$.
Then, \textsl{via} (\ref{FGdot1}), (\ref{basic12}) is equivalent to
\begin{align}
\label{zdotforcedyeah1} \dot{z}=\varepsilon
g_{-1}(\varepsilon)+\mu_1
 e^{ i t} H_1((z+F_G(t,\varepsilon)) e^{- i
t},\mbox{c.c.},\mu_1).
\end{align}
Generically, $g_{-1}(\varepsilon)\neq 0$ in some neighbourhood of
the origin. Thus, (\ref{zdotforcedyeah1}) does not generically have
a periodic solution when $\mu_1=0$ and $\varepsilon\neq 0$.
Furthermore, the analysis of (\ref{zdotforcedyeah1}) reduces to that
of (\ref{basiceqs5}) when $\varepsilon=0$ (see section
\ref{thecasee=0} and \cite{LW}). Set
$\varepsilon=\hat{\varepsilon}\mu_1$. Then (\ref{zdotforcedyeah1})
transforms to
\begin{equation}
\dot{z}=\mu_1 H_*(z e^{- i t},\ov{z} e^{ i
t},t,\hat{\varepsilon},\mu_1), \label{newzdotforced1}
\end{equation}
where
$H_*(w,\ov{w},t,\hat{\varepsilon},\mu_1)=\hat{\varepsilon}g_{-1}(\hat{\varepsilon}\mu_1)+
e^{ i t} H_1(w+F_G(t,\hat{\varepsilon}\mu_1) e^{- i t},\ccc,\mu_1)$
is smooth, $2\pi-$periodic in $t$ and uniformly bounded in $w$. \par
Set $\alpha_1=D_1H(- i v, i\ov{v},0,0)$. Near
$(z,\hat{\varepsilon},\mu_1)=(0,0,0)$, the time$-2\pi$ map $P$ of
(\ref{newzdotforced1}) is given by
\begin{equation}\label{Poincare234}
P(z,\hat{\varepsilon},\mu_1)=z+2\pi\mu_1\Big[\hat{\varepsilon}g_{-1}(0)+\alpha_1z+
O\big(|z|^2\big)+ O\big(\mu_1\big)+\mbox{higher order terms}\Big].
\end{equation} The fixed points of (\ref{Poincare234}) correspond to $2\pi-$periodic
solutions of (\ref{newzdotforced1}), and so to discrete rotating waves with
trivial spatio-temporal symmetry in (\ref{basic12}). Let
$$B_*(z,\overline{z},\hat{\varepsilon},\mu_1)=\hat{\varepsilon}g_{-1}(0)+\alpha_1z+ O\big(|z|^2\big)+ O(\mu_1)+\mbox{higher
order terms}$$ be the function inside the square brackets in
(\ref{Poincare234}). Note that $B_*(0,0,0,0)=0$ and that
$D_1B_*(0,0,0,0)=\alpha_1$, which is generically not $0$. By the
implicit function theorem, there is a unique smooth function
$z(\hat{\varepsilon},\mu_1)$ defined near
$(\hat{\varepsilon},\mu_1)=(0,0)$ with $z(0,0)=0$ and
\begin{align*}B_*\big(z(\hat{\varepsilon},\mu_1),\overline{z}(\hat{\varepsilon},\mu_1),\hat{\varepsilon},\mu_1\big)\equiv
0\end{align*} near $z=0$. However $z=0$ is not in general a solution
of $B_*=0$ for small $(\hat{\varepsilon},\mu_1)$, unless $g_{-1}(0)=
0$. Thus the generic behaviour of discrete rotating waves in
(\ref{basic12}) is to locally drift away from the origin. This leads
to the following theorem.
 \begin{thm} \label{noclueagain}Let $\alpha_1$ and $g_{-1}(0)\neq 0$ be as in the preceding discussion,
   with $\real(\alpha_1)\neq 0.$ Then there exists a wedge-shaped region near
    $(\varepsilon,\mu_1)=(0,0)$ of the form
$$
{\cal W}=\{(\varepsilon,\mu_1)\in
\GR^2\,:\,|\varepsilon|<K|\mu_1|,\,\,\,K>0,\,\,\mbox{$\mu_1$
near}\,\,0\,\}
$$
such that for all $(\varepsilon,\mu_1)\in {\cal W}$ with
$\varepsilon\neq 0$, $(\ref{basic12})$ has a unique hyperbolic
discrete rotating wave ${\cal D}^{1}_{\varepsilon,\mu_1}$ with
trivial spatio-temporal symmetry and center $[{\cal
D}^{1}_{\varepsilon,\mu_1}]_{\A}$ near, but generically not at, the
origin. Furthermore, $[{\cal D}^{1}_{\varepsilon,\mu_1}]_{\A}$ is a
center of anchoring when $\mu_1\real(a_1)< 0$.
\end{thm}
 \noindent\proof Let $z(\hat{\varepsilon},\mu_1)$ be the unique continuous function
solving the equation $B_*=0$ for small parameter vectors
$(\hat{\varepsilon},\mu_1)$, as asserted above. When $\mu_1=0$, any
$z$ is a non-hyperbolic fixed point of $P$ and so, from now on, we
will assume that $\mu_1\neq 0$. If that is the case, and if
$\hat{\varepsilon}$ and $\mu_1$ are small enough, the eigenvalues
$\omega_{1,2}(\hat{\varepsilon},\mu_1)$ of
$DP(z(\hat{\varepsilon},\mu_1),\hat{\varepsilon},\mu_1)$ satisfy
\begin{align*}
\left|\omega_{1,2}(\hat{\varepsilon},\mu_1)
\right|^2&=1+4\pi\mu_1\real(\alpha_1)+\mu_1
O(\hat{\varepsilon},\mu_1)\neq 1,
\end{align*} since $\real(\alpha_1)\neq 0$. When
$\mu_1\real(\alpha_1)<0$, the eigenvalues lie inside the unit circle
and the fixed point is asymptotically stable; otherwise, it is
unstable.
\par Let $K> 0$ be such that the preceding discussion holds whenever $\hat{\varepsilon}\in (-K,K)$ and let
$\mathcal{W}$ be as stated in the hypothesis. If
$(\hat{\varepsilon},\mu_1)$ is such that the time$-2\pi$ map
(\ref{Poincare234}) has a hyperbolic fixed point
$z(\hat{\varepsilon},\mu_1)$ near $0$, then (\ref{basic12}) has a
hyperbolic $2\pi-$periodic orbit
$\tilde{z}_{\hat{\varepsilon},\mu_1}(t)$ centered near the origin.
\par Let $0\neq\hat{\varepsilon}\in (-K,K)$  and set
$\varepsilon=\hat{\varepsilon}\mu_1$. Then $(\varepsilon,\mu_1)\in
\mathcal{W}$, as
$|\varepsilon|=|\hat{\varepsilon}||\mu_1|<K|\mu_1|,$ and
$\tilde{z}_{\hat{\varepsilon},\mu_1}(t)$ is a $2\pi-$periodic orbit
for the pair $(\varepsilon,\mu_1)$, which we denote by
$z_{\varepsilon,\mu_1}(t)$. Since $p=z-F_G$, (\ref{basic12}) has a
unique perturbed rotating wave $\mathcal{D}_{\varepsilon,\mu_1}^1$,
with
$$[\mathcal{D}_{\varepsilon,\mu_1}^1]_{\A}=\frac{1}{2\pi}\int_{0}^{2\pi}\!\!\!\!\left(z_{\varepsilon,\mu_1}(t)-F_G(t,\varepsilon)\right)\,
 d t=[z_{\varepsilon,\mu_1}]_{\A}.$$ But
$[z_{\varepsilon,\mu_1}]_{\A}= O(\varepsilon,\mu_1)$ as
$(\varepsilon,\mu_1)\to 0$ and so
$\mathcal{D}_{\varepsilon,\mu_1}^1$ is a discrete rotating wave with
trivial spatio-temporal symmetry and center near (but generically
not at) the origin. The conclusion about the stability of
$\mathcal{D}_{\varepsilon,\mu_1}^1$ follows directly from the
hyperbolicity of the eigenvalues of $DP$. \qed\newl When the
parameter values stray outside of ${\cal W}$, all that can
generically be said with certainty is that solutions of
(\ref{basic12}) locally drift away from the origin, and so anchoring
cannot be centered there. After drifting, the spiral may very well
get anchored at some point far from the origin, depending on the
global nature of the function $H_1$ in (\ref{basic12}).\par The
analysis of the situation near the $\varepsilon-$axis involves fixed
points of (\ref{Poincare234}) at $\infty$, which correspond to
traveling waves in (\ref{basic12}) \cite{LW,LeB}. Such an analysis
lies outside the scope of this paper; the situation will be
investigated at a later date.
\subsection{The Case $\jmath^*>1$}
Let $C^{0}_{\GR}(\GC)$ and $C^{1}_{\GR}(\GC)$ be the spaces of continuous and
continuously differentiable functions from $\GR$ to $\GC$, respectively, and
$\mathfrak{P}^{2\pi/\jmath^*}_t$ be the space of $2\pi/\jmath^*-$periodic
functions of the variable $t$. \par Then $C^0_{2\pi/\jmath^*}=\{f: f\in
\mathfrak{P}^{2\pi/\jmath^*}_t\cap C^0_{\GR}(\GC)\}$ and
$C^1_{2\pi/\jmath^*}=\{f:f\in \mathfrak{P}^{2\pi/\jmath^*}_t\cap
C^1_{\GR}(\GC)\}$ are Banach spaces when endowed with the respective norms
$$||u||_0=\sup\{|u(x)|:x\in [0,\textstyle{2\pi/\jmath^*}]\}\quad
\mbox{and}\quad ||u||_1=||u||_0+||u'||_0,$$ and the linear operator
${\cal Y}:C^1_{2\pi/\jmath^*}\to C^0_{2\pi/\jmath^*}$ defined by
${\cal Y}(u)=iu+u'$ is bounded, invertible and has bounded inverse.
 The nonlinear operator ${\cal
H}_G:C^1_{2\pi/\jmath^*}\times {\mathbb R}^2\to C^0_{2\pi/\jmath^*}$
given by
\begin{equation}\label{HG}
{\cal H}_G(u,\varepsilon,\mu_1)={\cal Y}(u)-\mu_1 H_1\left(u- i v+
\varepsilon\sum_{m\in {\mathbb Z}}\,\frac{g_m(\varepsilon) e^{ i
m\jmath^* t}}{ i(m\jmath^*+1)}, \mbox{c.c.},\mu_1\right),
\end{equation}
where the $g_m$ are the Fourier coefficients of $G$ found in
(\ref{FourierG1}), will play an important part in what follows. Note
that ${\cal H}_G(0,0,0)=0$ and $D_1{\cal H}_G(0,0,0)= i\neq 0$.
Thus, by the implicit function theorem, there is a neighbourhood
${\cal N}$ of the origin in ${\mathbb R}^2$ and a unique smooth
function $U:{\mathbb R}^2\to C^1_{2\pi/\jmath^*}$ satisfying
$U(0,0)=0$ and ${\cal
H}_G(U(\varepsilon,\mu_1),\varepsilon,\mu_1)\equiv 0$ for all
$(\varepsilon,\mu_1)\in {\cal N}$.\newl Let $F_G:\GR\times\GR^2\to
\GC$ be defined by
\begin{equation}\label{FG}
F_G(t,\varepsilon,\mu_1)= e^{ i t}\Big[- i v+\varepsilon \sum_{m\in
{\mathbb Z}}\,\frac{g_m(\varepsilon) e^{ i m\jmath^* t}}{
i(m\jmath^*+1)}+ U(\varepsilon,\mu_1)(t)\Big].
\end{equation}
Then
\begin{align}\label{YU}
{\cal Y}(U(\varepsilon,\mu_1))(t)=\mu_1
H_1\left(F_G(t,\varepsilon,\mu_1) e^{- i t},\mbox{c.c.},\mu_1\right)
\end{align}
and
\begin{align}\label{FGdot}
\dot{F}_G(t,\varepsilon,\mu_1)&= e^{ i t}\big[v+\varepsilon
G(t,\varepsilon,\mu_1)+{\cal
Y}\left(U(\varepsilon,\mu_1)\right)(t)\big].
\end{align}
On the center bundle $\GC\times GS^1$, the path
$(F_G(t,\varepsilon,\mu_1),t)$ (generically) represents a discrete
rotating wave around the origin with
$\GZ_{\jmath^*}-$spatio-temporal symmetry. Set $z=p-F_G$. Then,
using (\ref{YU}) and (\ref{FGdot}), (\ref{basic12}) rewrites as
\begin{align}
\label{zdotforcedyeah} \dot{z}=\mu_1  e^{ i t}\big[
H_1((z+F_G(t,\varepsilon,\mu_1)) e^{- i
t},\mbox{c.c.},\mu_1)-H_1(F_G(t,\varepsilon,\mu_1) e^{- i
t},\mbox{\mbox{c.c.}},\mu_1))\big],
\end{align}
which reduces to
\begin{equation}
\dot{z}=\mu_1  e^{ i t}\widehat{H}(z e^{- i t},\ov{z} e^{ i
t},t,\varepsilon,\mu_1), \label{newzdotforced}
\end{equation}
where
$$\widehat{H}(w,\ov{w},t,\varepsilon,\mu_1)=H_1(w+F_G(t,\varepsilon,\mu_1) e^{- i t},\mbox{c.c.},\mu_1)-H_1(F_G(t,\varepsilon,\mu_1) e^{- i t},\mbox{c.c.},\mu_1)$$
is $2\pi/\jmath^*-$periodic in $t$.\par Near
$(z,\varepsilon,\mu_1)=(0,0,0)$, the time$-2\pi$ map $P$ of
(\ref{newzdotforced}) is given by
\begin{equation}\label{Poincare23}
P(z,\varepsilon,\mu_1)=z+2\pi\mu_1\Big[(\alpha_1+
O(\varepsilon,\mu_1))z+ O(\varepsilon,\mu_1)\ov{z}+
O\big(|z|^2\big)\Big],
\end{equation} where $\alpha_1$ is as defined in section~\ref{thecasel=1}. The fixed points of (\ref{Poincare23}) correspond to $2\pi-$periodic
solutions of (\ref{newzdotforced}), and so to discrete rotating waves with
$\GZ_{\jmath^*}-$spatio-temporal symmetry in (\ref{basic12}). Let
\begin{equation*}B^*(z,\overline{z},\varepsilon,\mu_1)=(\alpha_1+ O(\varepsilon,\mu_1))z+ O(\varepsilon,\mu_1)\ov{z}+ O\big(|z|^2\big)\end{equation*} be the function inside the square brackets in (\ref{Poincare23}). Note that
$B^*(0,0,0,0)=0$ and that $D_1B^*(0,0,0,0)=\alpha_1$, which is
generically not $0$. By the implicit function theorem, there is a
unique smooth function $z(\varepsilon,\mu_1)$ defined near
$(\varepsilon,\mu_1)=(0,0)$ with $z(0,0)=0$ and
\begin{align*}B^*\big(z(\varepsilon,\mu_1),\overline{z}(\varepsilon,\mu_1),\varepsilon,\mu_1\big)\equiv
0\end{align*} near $z=0$. But $z=0$ is always a solution of $B^*=0$;
hence $z(\varepsilon,\mu_1)=0$ for all small enough
$(\varepsilon,\mu_1)$. This leads to the following theorem.
 \begin{thm} \label{noclue}Let $\alpha_1$ be as in the preceding discussion,
   with $\mbox{\rm Re}(\alpha_1)\neq 0.$ If $(\varepsilon,\mu_1)$ is in a small deleted neighbourhood
   $\mathcal{W}^{\jmath^*}$ of the origin, $(\ref{basic12})$ has a unique hyperbolic
discrete rotating wave ${\cal D}^{\jmath^*}_{\varepsilon,\mu_1}$
with $\GZ_{\jmath^*}-$spatio-temporal symmetry and center $[{\cal
D}^{\jmath^*}_{\varepsilon,\mu_1}]_A=0$. Furthermore, ${\cal
D}^{\jmath^*}_{\varepsilon,\mu_1}$ is anchored at the origin if
$\real(\alpha_1)\mu_1<0$.
\end{thm}


\section{Spiral Anchoring With General FESB}\label{general}

The case $n>1$ combines the results from the preceding section and from section
3.2 in \cite{BLM}; it provides a synopsis of spiral anchoring under the most
general form of Euclidean symmetry-breaking. The results depend yet again on
the nature of $\jmath^*$. The proofs follow the general lines of previous work
and are thus omitted, in order to avoid unnecessary repetitions.

\begin{thm} \label{noclueagainconj}Let ${\jmath^*}=1$  and  $k\in
\{1,\ldots,n\}$. Set $\alpha_k=D_1H_k(- i v, i \ov{v},0,0)$. Write
$\varepsilon=\mu_0$. Generically, there is a wedge-shaped region of
the form
$$
{\mathcal W}^1_k=\{(\mu_0,\mu_1,\ldots,\mu_n)\in
\GR^{n+1}\,:\,|\mu_j|<W_{j,k}|\mu_k|,\,\,\,W_{j,k}>0,\,\,\mbox{\rm for $j\neq
k$ and $\mu_k$ near}\,\,0\,\}
$$
such that for all $(\varepsilon,\mu_1,\ldots,\mu_n)\in {\mathcal W}_{k}^1$ with
$\varepsilon\neq 0$, the center bundle equations $(\ref{basic12})$ have a
unique hyperbolic discrete rotating wave
$\mathfrak{D}^{1;k}_{\varepsilon,\mu}$, with trivial spatio-temporal symmetry,
centered away from $\xi_k$. Furthermore,
$[\mathfrak{D}^{1;k}_{\varepsilon,\mu}]_{\A}$ is a center of anchoring when
$\mu_k\real(\alpha_k)<0$.
\end{thm}
\begin{thm} \label{noclueconj}Let ${\jmath^*}>1$ and $k\in \{1,\ldots,n\}$. Set $\alpha_k=D_1H_k(- i v, i \ov{v},0,0)$.
Ge\-ne\-ri\-cal\-ly, there is a cone-like region of the form \small
$$
{\mathcal W}^{\jmath^*}_k=\{(\varepsilon,\mu_1,\ldots,\mu_n)\in
\GR^{n+1}\,:\,|\varepsilon|<\varepsilon_0,\,\,\,|\mu_j|<W_{j,k}|\mu_k|,\,\,\,W_{j,k}>0,\,\,\mbox{\rm
for $j\neq k$ and $\mu_k$ near}\,\,0\,\}
$$\normalsize
such that for all $(\varepsilon,\mu)\in {\mathcal W}^{\jmath^*}_k$,
$(\ref{basic12})$ has a unique hyperbolic discrete rotating wave
 $\mathfrak{D}^{{\jmath^*};k}_{\varepsilon,\mu}$ with $\GZ_{{\jmath^*}}-$spatio-temporal symmetry centered at
 $\xi_k$. Furthermore,
$[\mathfrak{D}^{\jmath^*;k}_{\varepsilon,\mu}]_{\A}$ is a center of anchoring
when $\mu_k\real(\alpha_k)<0$.
\end{thm}


\section{Numerical Simulations}
In this section, we examine a system that models excitable media
with anisotropy and a single inhomogeneity, which give rise to a
semi-flow satisfying the FESB equivariance conditions described by
(E1), (E2'), (E3) and (E4) for $\jmath^*=1$.
\par The computations are carried out on a two-dimensional square
domain\footnote{More precisely, on the domain $[-30,30]^2$ with 200
grid points to a side and time-step $\Delta t=0.005$.} with Neumann
boundary condition, using a 5-point Laplacian and the Runge-Kutta
forward integrator of order 2. This
 naive numerical approach is used instead of a more robust algorithm mostly because the aim of this section is to
illustrate the qualitative (rather than quantitative) power of theorems
\ref{noclueagainconj} and \ref{noclueconj}.
\par Given appropriate initial conditions, this system sustains
spiral waves. A planar reaction-diffusion-advection system (RDAS) is
a reaction-diffusion system to which an advection term has been
added:
\begin{equation}\label{rdacmrt} u_t=\tilde{D}\Delta
u+\tilde{A}_1u_{x_1}+\tilde{A}_2u_{x_2}+f(u),\end{equation} where
$\tilde{A}_i\in \GM_2(\GR)$ for $i=1,2$, and the other terms are as in
(\ref{basicrdpde}). \par Advection terms are used to model a directed flow or
current through the excitable medium under consideration;
 for instance, the modified RDAS \begin{equation}
\label{rdapert} u_t=\tilde{D}\Delta
u+f(u)+\varepsilon\left[\tilde{A}_1(\varepsilon,\mu)u_{x_1}+\tilde{A}_2(\varepsilon,\mu)u_{x_2}\right]+\mu
g(u,\|x-\xi\|^2,\varepsilon,\mu),
\end{equation} where $(\varepsilon,\mu)\in \GR^2$ is small, $\tilde{A}_1,\tilde{A}_2:\GR^2\to \GM_2(\GR)$ are smooth functions and $g$ is a smooth bounded function, could model membrane potentials in a piece of cardiac tissue subjected to a directed current, with an inhomogeneity located at $\xi$. \par Let $\Phi_{t,\varepsilon,\mu}$
denote the semi-flow generated by (\ref{rdapert}): it is
$\GSE(2)-$equi\-va\-ri\-ant when $(\varepsilon,\mu)=0$; $\GC\sdp
\{1\}-$equi\-va\-riant when $\varepsilon\neq 0$ is small and
$\mu=0$; $\GSO(2)_{\xi}-$equi\-va\-ri\-ant when $\mu\neq 0$ is small
and $\varepsilon=0$, and (generically) trivially equivariant when
$(\varepsilon,\mu)$ is small and generic.\footnote{This can be seen
by slightly modifying the proof of theorem 2.1.3 in
\cite{Byeah}.}\par However, it is not clear if
theorem~\ref{noclueagain} can be applied to (\ref{rdapert}) since it
is not yet known if the center manifold reduction theorems of
\cite{SSW1,SSW2,SSW3,SSW4} hold in general for RDAS; as such, the
center bundle equations need not be of the form (\ref{basic10}) in
the presence of advection.\par Yet, the following example shows that
the conclusion of theorem~\ref{noclueagain} holds for the following
RDAS, and so that, in some sense, (\ref{basic10}) may capture the
essential dynamics near a rotating wave when the semi-flow of the
PDE has the general symmetry-breaking properties outlined above.
\newl Consider the RDAS
\begin{align}\label{FHNM2RDAS}
\begin{split}
  u_t&= \textstyle{\frac{1}{\varsigma}
  \left(u-\frac{1}{3}u^3-v\right)-3\sqrt{2}\sin\left(\frac{0.03\pi}{2}\right)\phi_2+\Delta u+0.002u_{x_1}},\\
  v_t&=\varsigma(u+\beta-\gamma v+\phi_2),
\end{split}
\end{align} where $\varsigma=0.3$, $\beta=0.6$, $\gamma=0.5$ and
$\phi$ is defined \textit{via}
\begin{align*}\phi(x) =0.12 f(x_1+10,x_2-5\sqrt{3}),\end{align*} with $f(x)=\exp\left(-0.00086\left(x_1^2+x_2^2\right)\right).$
 A single integration (see
\begin{figure}[t]\begin{center}
\includegraphics[height=175pt]{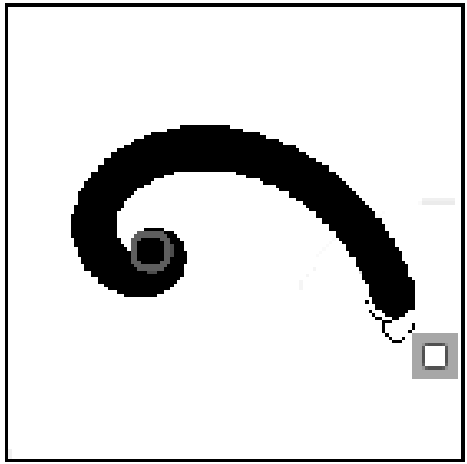}
\caption{Anchoring in the RDAS equations (\ref{FHNM2RDAS}). The
spiral tip path is plotted in black. The anchored perturbed rotating
wave is shown as a gray closed curve; the gray square indicates the
location of the perturbation center.}
\end{center}\hrule\end{figure}\ figure~\thefigure) clearly shows the spiral anchoring away from the perturbation center.

\section*{References}
\bibliographystyle{amsplain}
\bibliography{spiralsp}
\end{document}